\documentclass[reqno]{amsart}
\usepackage{amsmath}
\usepackage{amssymb}
\usepackage{amsthm}
\usepackage{amsfonts}
\usepackage{mathtools}
\usepackage{bm}
\usepackage{ulem}
\usepackage{fancyhdr}
\usepackage{hyperref}
\usepackage[pdftex,final]{graphicx}
\usepackage[iso-8859-7]{inputenc}

\theoremstyle{plain}

\newtheorem*{theorem*}{Theorem}
\theoremstyle{definition}

\numberwithin{equation}{section}

\newcommand{\mytilde}{\raise.17ex\hbox{$\scriptstyle\mathtt{\sim}$}}

\newcommand{\RgeO}{\ensuremath{\mathbb{R}_{ \geq 0}}}
\newcommand{\Rat}[1]{\ensuremath{\mathbb{R}^{#1}}}

\begin{document}
\title{Observer design for a general class of triangular systems}

\author{D. Boskos}
\address{Department of Mathematical and Physical Sciences,
        National Technical University of Athens, Zografou Campus 15780, Athens, Greece}
\email{dmposkos@central.ntua.gr}
%\ead{ykarag@math.ntua.gr} 
 %\author{N. Yannakakis\corref{cor1}}%\fnref{label1}

%\ead{nyian@math.ntua.gr} 

\author{J. Tsinias}

%\cortext[cor1]{Corresponding author}
%\fntext[label1]{blaa}
\address{Department of Mathematical and Physical Sciences,
        National Technical University of Athens, Zografou Campus 15780, Athens, Greece}
\email{jtsin@central.ntua.gr}

\begin{abstract}
The paper deals with the observer design problem for a wide class of triangular nonlinear systems. Our main results generalize those obtained in the recent author's works \cite{BdTj13a} and \cite{BdTj13b}.     
\end{abstract}

%% MSC codes here, in the form: 
%\subjclass[2000]{ 35J20, 35J60, 35J70}.
%% or \MSC[2008] code \sep code (2000 is the default)

\maketitle
\section{INTRODUCTION}

%\newcounter{numb}
%\setcounter{numb}{1}

We derive sufficient conditions for the solvability of the observer design problem for time-varying single output systems of the form
\begin{subequations} \label{system:triangular}
\begin{align}
\dot{x}_{i}&=f_{i}(t,x_{1} ,...,x_{i+1}),i=1,...,n-1 \nonumber \\
\dot{x}_{n}&=f_{n} (t,x_{1} ,...,x_{n}), \label{system:triangular:st} \\
y&=x_{1}, (x_{1} ,...,x_{n})\in\Rat{n} \label{system:triangular:out}
\end{align} 
\end{subequations}

\noindent where the functions $f_{i} (\cdot )$, $i=1,2,...,n$, are continuous and (locally) Lipschitz. It is known from \cite{GjKi01} that every single output control system which has a uniform canonical flag (\cite[Chapter 2, Definition 2.1]{GjKi01}) can be locally transformed in the above canonical form \eqref{system:triangular} for each fixed input. We also mention the works \cite{CdLw03} and \cite{Rw03} where algebraic type necessary and sufficient conditions are established for feedback equivalence between a single input system and a triangular system whose dynamics have $p$-normal form. In our recent work \cite{BdTj13b}, the observer design problem is studied for a subclass of systems \eqref{system:triangular} whose dynamics have $p$-normal form. The result of present work constitutes a generalization of previous results in the literature dealing with the observer design problem for triangular systems (see for instance \cite{AvPlAa09}, \cite{GjKi01}, \cite{HhTbAf02} and relative references therein) and particularly generalizes the main result of the recent author's work in \cite{BdTj13b}. % , \cite{SrPl11}

We make the following assumption for the right hand side of system \eqref{system:triangular}.

\textbf{H1.} For each $(t;x_{1},...,x_{i})\in {\mathbb R}_{\ge 0}\times {\mathbb R}^{i}$, $i=1,...,n-1$, the function $\mathbb{R}\ni z\to f_{i}(t,x_{1} ,...,x_{i},z)\in\mathbb{R}$ is strictly monotone.

The paper is organized as follows. We first provide notations and various concepts, including the concept of the \textit{switching observer} that has been originally introduced in \cite{BdTj13a}, for general time-varying systems:
\begin{subequations} \label{system:nonlinear}
\begin{align}
\dot{x}&=f(t,x),(t,x)\in {\mathbb R}_{\ge 0} \times {\mathbb R}^{n} \label{system:nonlinear:st} \\
y&=h(t,x),y\in {\mathbb R}^{\bar{n}} \label{system:nonlinear:out}
\end{align} 
\end{subequations}

\noindent where $y(\cdot )$ plays the role of output. We then provide the precise statement of our main result (Theorem 1.1) concerning solvability of the observer design problem for \eqref{system:triangular}. Section II contains some preliminary results concerning solvability of the observer design problem for the case \eqref{system:nonlinear} with linear output (Propositions 2.1 and 2.2). In Section III we use the results of Section II, in order to prove our main result.  

\textit{Notations and definitions:} We adopt the following notations. For a given vector $x\in{\mathbb R}^{n}$, $x'$ denotes its transpose and $|x|$ its Euclidean norm. We use the notation $|A|:=\max \{|Ax|:x\in {\mathbb R}^{n};|x|=1\}$ for the induced norm of a matrix $A\in {\mathbb R}^{m\times n} $. By $N$ we denote the class of all increasing $C^{0} $ functions $\phi:{\mathbb R}_{\ge 0}\to{\mathbb R}_{\ge 0} $. For given $R>0$, denote by $B_{R}$ the closed ball of radius $R>0$, centered at $0\in{\mathbb R}^{n} $. Consider a pair of metric spaces $X_{1}$, $X_{2}$ and a set-valued map $X_{1}\ni x\to Q(x)\subset X_{2}$. We say that $Q(\cdot)$ satisfies the \textit{Compactness Property }(\textbf{CP)}, if for every sequence $(x_{\nu})_{\nu \in {\mathbb N}}\subset X_{1}$ and $(q_{\nu})_{\nu \in {\mathbb N}}\subset X_{2}$ with $x_{\nu }\to x\in X_{1}$ and $q_{\nu} \in Q(x_{\nu})$, there exist a subsequence $(x_{\nu_{k}} )_{k\in{\mathbb N}}$ and $q\in Q(x)$ such that $q_{\nu_{k}}\to q$. We also invoke the well known fact, see \cite{Ki05}, that the time-varying system \eqref{system:nonlinear:st} is forward complete, if and only if there exists a function $\beta \in NN$ such that the solution $x(\cdot ):=x(\cdot ,t_{0} ,x_{0} )$ of \eqref{system:nonlinear:st} initiated from $x_{0} $ at time $t=t_{0} $ satisfies:
\begin{equation} \label{state:bound} 
|x(t)|\le\beta (t,|x_{0}|), \forall t\ge t_{0}\ge 0, x_{0} \in {\mathbb R}^{n}  
\end{equation}

\noindent provided that the dynamics of \eqref{system:nonlinear:st} are $C^{0}$ and  Lipschitz on $x\in {\mathbb R}^{n} $. It turns out that, under these regularity assumptions plus forward completeness for \eqref{system:nonlinear:st}, for each $t_{0}\ge 0$ and $x_{0}\in {\mathbb R}^{n}$ the corresponding output $y(t)=h(t,x(t,t_{0},x_{0}))$ of \eqref{system:nonlinear} is defined for all $t\ge t_{0}$. For each $t_{0} \ge 0$ and nonempty subset $M$ of ${\mathbb R}^{n} $, we may consider the set $O(t_{0},M)$ of all outputs of \eqref{system:nonlinear}, corresponding to initial state $x_{0} \in M$ and initial time $t_{0} \ge 0$:
\begin{align} \label{output:functions} 
O(t_{0},M):=\{y&:[t_{0},\infty)\to {\mathbb R}^{\bar{n}} \nonumber \\
:y&(t)=h(t,x(t,t_{0},x_{0}));t\ge t_{0},x_{0} \in M\}  
\end{align} 

\noindent For given \textit{$\emptyset \ne M\subset {\mathbb R}^{n}$}, we say that the \textbf{Observer Design Problem (ODP) is solvable for \eqref{system:nonlinear} with respect to $M$}, if for every $t_{0} \ge 0$ there exist a continuous map $G:=G_{t_{0}}(t,z,w):[t_{0},\infty)\times {\mathbb R}^{n}\times {\mathbb R}^{\bar{n}}\to {\mathbb R}^{n}$ and a nonempty set $\bar{M}\subset {\mathbb R}^{n}$ such that for every $z_{0} \in \bar{M}$ and output $y\in O(t_{0},M)$ the corresponding trajectory $z(\cdot):=z(\cdot,t_{0},z_{0};y)$; $z(t_{0})=z_{0}$ of the observer $\dot{z}=G(t,z,y)$ exists for all $t\ge t_{0} $ and the error $e(\cdot):=x(\cdot)-z(\cdot)$ between the trajectory $x(\cdot):=x(\cdot,t_{0},x_{0})$, $x_{0} \in M$ of \eqref{system:nonlinear:st} and the trajectory $z(\cdot ):=z(\cdot,t_{0},z_{0};y)$ of the observer satisfies:
\begin{equation} \label{error:conv} 
\mathop{\lim }\limits_{t\to \infty } e(t)=0 
\end{equation}

\noindent We say that the \textbf{Switching Observer Design Problem (SODP) is solvable for \eqref{system:nonlinear} with respect to $M$}, if for every $t_{0} \ge 0$ there exist a strictly increasing sequence of times $(t_{k})_{k\in {\mathbb N}}$ with $t_{1}=t_{0}$ and $\lim_{k\to\infty} t_{k}=\infty$, a sequence of continuous mappings $G_{k} :=G_{k,t_{k-1}}(t,z,w):[t_{k-1},t_{k+1}]\times {\mathbb R}^{n}\times {\mathbb R}^{\bar{n}}\to {\mathbb R}^{n}$, $k\in {\mathbb N}$ and a nonempty set $\bar{M}\subset {\mathbb R}^{n}$ such that the solution $z_{k} (\cdot)$ of the system
\begin{equation} \label{observer:seq} 
\dot{z}_{k} =G_{k} (t,z_{k} ,y), t\in [t_{k-1} ,t_{k+1} ] 
\end{equation} 

\noindent with initial $z(t_{k-1} )\in \bar{M}$ and output $y\in {\rm O}(t_{0} ,M)$, is defined for every $t\in [t_{k-1} ,t_{k+1} ]$ and in such a way that, if we consider the piecewise continuous map $Z:[t_{0} ,\infty )\to {\mathbb R}^{n} $ defined as $Z(t):=z_{k}(t)$, $t\in [t_{k},t_{k+1})$, $k\in {\mathbb N}$, where for each $k\in {\mathbb N}$ the map $z_{k}(\cdot)$ denotes the solution of \eqref{observer:seq}, then the error $e(\cdot):=x(\cdot)-Z(\cdot)$ between the trajectory $x(\cdot ):=x(\cdot ,t_{0} ,x_{0} )$, of \eqref{system:nonlinear:st} and $Z(\cdot )$ satisfies \eqref{error:conv}.

\noindent Our main result is the following theorem.

\textit{Theorem 1.1:} For the system \eqref{system:triangular}, assume that Hypothesis H1 is satisfied and \eqref{system:triangular:st} is forward complete, i.e. there exists a function $\beta \in NN$, such that the solution $x(\cdot ):=x(\cdot ,t_{0} ,x_{0} )$ of \eqref{system:triangular:st} satisfies the estimation \eqref{state:bound}. Then:

\noindent (i) the SODP is solvable for \eqref{system:triangular} with respect to ${\mathbb R}^{n} $.

\noindent (ii) if in addition we assume that it is a priori known, that the initial states of \eqref{system:triangular} belong to a given ball $B_{R} $ of radius $R>0$ centered at zero $0\in {\mathbb R}^{n} $, then the ODP is solvable for \eqref{system:triangular} with respect to $B_{R} $. \ensuremath{\triangleleft}

Theorem 1.1 constitutes a generalization of Theorem 1.1 in \cite{BdTj13b} for systems \eqref{system:triangular}, whose dynamics are $C^{1}$ and have the particular form $f_{i}(t,x_{1},...,x_{i+1}):=\tilde{f}_{i}(t,x_{1},...,x_{i})+a_{i}(t,x_{1})x_{i+1}^{m_{i}}$, for certain $\tilde{f}_{i}\in C^{1}(\RgeO\times\Rat{i};\mathbb{R})$ and $a_{i}\in C^{1}(\RgeO\times\Rat{};\mathbb{R})$, $i=1,...,n-1$, where the constants $m_{i}$, $i=1,...,n-1$ are odd integers and the functions $a_{i}(\cdot,\cdot)$, $i=1,...,n-1$ satisfy the condition $|a_{i}(t,y)|>0$, $\forall t\in\RgeO$, $y\in\mathbb{R}$. 

\section{PRELIMINARY RESULTS}

The proof of our main result concerning the case \eqref{system:triangular}, is based on some preliminary results concerning the case of systems \eqref{system:nonlinear} with linear output:
\begin{subequations} \label{system:linearout}
\begin{align}
\dot{x}&=f(t,x):=F(t,x,H(t)x),(t,x)\in {\mathbb R}_{\ge 0}\times {\mathbb R}^{n} \label{system:linearout:st} \\
y&=h(t,x):=H(t)x,\; y\in {\mathbb R}^{\bar{n}} \label{system:linearout:out}
\end{align}
\end{subequations}

\noindent where $H:{\mathbb R}_{\ge 0}\to {\mathbb R}^{\bar{n}\times n}$ is $C^{0}$ and $F:{\mathbb R}_{\ge 0}\times {\mathbb R}^{n}\times {\mathbb R}^{\bar{n}}\to {\mathbb R}^{n}$ is $C^{0}$ and Lipschitz on $(x,y)\in {\mathbb R}^{n}\times {\mathbb R}^{\bar{n}} $. We assume that system \eqref{system:linearout:st} is forward complete, namely, the solution $x(\cdot ):=x(\cdot ,t_{0} ,x_{0} )$ of \eqref{system:linearout:st} satisfies \eqref{state:bound} for certain $\beta \in NN$, hence, for every $R>0$ and $t\ge 0$ we can define:
\begin{equation} \label{output:set} 
Y_{R}(t):=\{y\in {\mathbb R}^{\bar{n}} :y=H(t)x(t,t_{0} ,x_{0} ),\:{\rm for}\:{\rm certain}\: t_{0} \in [0,t]\:{\rm and}\: x_{0}\in B_{R}\} 
\end{equation}

\noindent where $H(\cdot )$ is given in \eqref{system:linearout:out}. Obviously, $Y_{R}(t)\ne \emptyset$ for all $t\ge 0$ and, if \eqref{state:bound} holds, the set-valued map $[0,\infty)\ni t\to Y_{R}(t) \subset {\mathbb R}^{\bar{n}}$ satisfies the CP and further $y(t)\in Y_{R}(t)$, for every $t\ge t_{0}\ge 0$ and $y\in O(t_{0},B_{R})$; (see notations). Also, given integers $\ell,m,n,\bar{n}\in\mathbb{N}$ and a map $A:{\mathbb R}_{\ge 0}\times {\mathbb R}^{\ell}\times {\mathbb R}^{n}\times {\mathbb R}^{m}\times {\mathbb R}^{\bar{n}}\to\Rat{m\times m}$, we say that $A(\cdot,\cdot,\cdot,\cdot,\cdot)$ satisfies Property P1, if the following holds:

\textbf{P1.} $A(t,q,x,e,y)$ has the form
\begin{equation}
A(t,q,x,e,y):=(A_{C1}(t,q,x,e_{1},y), A_{C2}(t,q,x,e_{1},e_{2},y), \cdots, A_{Cm}(t,q,x,e_{1},...,e_{m},y)) \label{map:A:aux:form} 
\end{equation}

\noindent where each mapping $A_{Ci}:{\mathbb R}_{\ge 0}\times {\mathbb R}^{\ell}\times {\mathbb R}^{n}\times \Rat{i}\times {\mathbb R}^{\bar{n}}\to\Rat{m\times 1}$, $i=1,...,m$ is continuous on ${\mathbb R}_{\ge 0}\times {\mathbb R}^{\ell}\times {\mathbb R}^{n}\times \{(e_{1},...,e_{i})\in\Rat{i}:e_{i}\ne 0\}\times {\mathbb R}^{\bar{n}}$  and bounded on every compact subset of ${\mathbb R}_{\ge 0}\times {\mathbb R}^{\ell}\times {\mathbb R}^{n}\times\Rat{i}\times {\mathbb R}^{\bar{n}}$.

\noindent We make the following hypothesis:

\textit{Hypothesis 2.1.} There exist a function $g\in C^{1}({\mathbb R}_{\ge 0};{\mathbb R})$ satisfying:
\begin{subequations} \label{map:g:properties}
\begin{align}
0<g(t)<1,&\forall t\ge 0; \label{map:g:values} \\
\dot{g}(t)\ge -g(t),&\forall t\ge 0; \\
\lim_{t\to \infty }g(t)&=0 \label{map:g:convergence}
\end{align}
\end{subequations}

\noindent an integer $\ell \in {\mathbb N}$, a map $A:{\mathbb R}_{\ge 0}\times {\mathbb R}^{\ell}\times {\mathbb R}^{n}\times {\mathbb R}^{n}\times {\mathbb R}^{\bar{n}}\to{\mathbb R}^{n\times n}$ satisfying P1 and constants $L>1$, $c_{1},c_{2}>0$, $R>0$ with
\begin{equation} \label{constant:c1} 
c_{1}\ge\tfrac{1}{2}  
\end{equation} 

\noindent such that the following properties hold:

\textbf{A1.} For every $\xi\ge 1$ there exists a set-valued map
\begin{equation} \label{map:Q} 
[0,\infty) \ni t\to Q_{R}(t):=Q_{R,\xi}(t)\subset {\mathbb R}^{\ell}  
\end{equation} 

\noindent with $Q_{R}(t)\ne \emptyset$ for any $t\ge 0$, satisfying the CP and such that
\begin{align}
&\forall t\ge 0,y\in Y_{R}(t),x,z\in {\mathbb R}^{n}\:{\rm with}\: |x|\le\beta (t,R)\: {\rm and}\: |x-z|\le\xi \nonumber \\
\Rightarrow &\Delta F(t,x,z,y):=F(t,x,y)-F(t,z,y)=A(t,q,x,x-z,y)(x-z)\: {\rm for}\: {\rm certain}\: q\in Q_{R}(t) \label{rhs:difference}                                 
\end{align}

\noindent with $Y_{R}(\cdot)$ as given by \eqref{output:set}.

\textbf{A2.} For every $\xi \ge 1$ there exists a set-valued map $Q_{R}:=Q_{R,\xi}$ as in Hypothesis A1, in such a way that for every $t_{0} \ge 0$, a time-varying symmetric matrix $P_{R}:=P_{R,\xi,t_{0}}\in C^{1}([t_{0},\infty);{\mathbb R}^{n\times n})$ and a function $d_{R}:=d_{R,\xi,t_{0}}\in C^{0}([t_{0},\infty);{\mathbb R})$ can be found, satisfying:

\begin{subequations} \label{A2:hypothesis}
\begin{equation}\label{map:PR:properties}
P_{R}(t)\ge I_{n\times n},\forall t\ge t_{0};|P_{R}(t_{0})|\le L;                                         
\end{equation}
\begin{equation} \label{map:dR:properties}
d_{R}(t)>c_{1}, \forall t\ge t_{0} +1; \int_{t_{1}}^{t_{2}}d_{R}(s)ds>-c_{2},\forall t_{2} \ge t_{1} \ge t_{0};                       
\end{equation}
\begin{align}
e'P_{R}(t)A(t,q,x,e,y)e+\tfrac{1}{2}e'\dot{P}_{R}(t)e\le -d_{R}(t)e'P_{R}(t)e,\forall t\ge t_{0},q&\in Q_{R}(t),\nonumber \\
x\in {\mathbb R}^{n},e\in\ker H(t),y\in Y_{R}(t):|x|\le\beta(t,R),|e|\le\xi,e'P_{R}(t)e&\ge g(t) \label{Lyapunov:ker:inequality}
\end{align}
\end{subequations}

\noindent with $Y_{R}(\cdot)$ as given by \eqref{output:set}.

The following result, constitutes a slight generalization of Proposition 2.1 in \cite{BdTj13b}.

\textit{Proposition 2.1:} Consider the system \eqref{system:linearout} and assume that it is forward complete, namely, \eqref{state:bound} holds for certain $\beta \in NN$ and satisfies Hypothesis 2.1. Then, the following hold:

\noindent For each $\bar{t}_{0} \ge t_{0} \ge 0$ and constant $\xi$ satisfying
\begin{equation}
\xi\ge\sqrt{L}\exp\{2c_{2}\}\beta(\bar{t}_{0},\bar{R}),\bar{R}:=R+1
\end{equation}

\noindent there exists a function $\phi_{R}:=\phi_{R,\xi,\bar{t}_{0}}\in C^{1}([\bar{t}_{0},\infty);{\mathbb R}_{>0})$, such that the solution $z(\cdot)$ of system
\begin{subequations} \label{observer:equation}
\begin{align}
\dot{z}=G_{\bar{t}_{0}}(t,z,y):=&F(t,z,y)+\phi_{R}(t)P_{R}^{-1}(t)H'(t)(y-H(t)z) \\
&{\rm with}\: {\rm initial}\: z(\bar{t}_{0})=0
\end{align}
\end{subequations}

\noindent where $P_{R}(\cdot ):=P_{R,\xi,\bar{t}_{0}}(\cdot)$ is given in A2, is defined for all $t\ge\bar{t}_{0}$ and the error $e(\cdot):=x(\cdot)-z(\cdot)$ between the trajectory $x(\cdot):=x(\cdot,t_{0},x_{0})$ of \eqref{system:linearout:st}, initiated from $x_{0}\in B_{R}$ at time $t_{0}\ge 0$ and the trajectory $z(\cdot):=z(\cdot,\bar{t}_{0},0;y)$ of \eqref{observer:equation} satisfies:
\begin{subequations}
\begin{align}
&|e(t)|<\xi,\forall t\ge\bar{t}_{0}; \label{error:bound} \\
&|e(t)|\le\max\{\xi\exp\{-c_{1}(t-(\bar{t}_{0}+1))\},\sqrt{g(t)}\},\forall t\ge \bar{t}_{0}+1 \label{error:ertimate}
\end{align}
\end{subequations}

\noindent It follows from \eqref{map:g:convergence} and \eqref{error:ertimate}, that for $\bar{t}_{0}:=t_{0}$ the ODP is solvable for \eqref{system:linearout} with respect to $B_{R}$; particularly the error $e(\cdot)$ between the trajectory $x(\cdot):=x(\cdot,t_{0},x_{0})$, $x_{0}\in B_{R}$ of \eqref{system:linearout:st} and the trajectory $z(\cdot ):=z(\cdot ,t_{0} ,z_{0} ;y)$, $z_{0} =0$ of the observer $\dot{z}=G_{t_{0} } (t,z,y)$ satisfies \eqref{error:conv}. \ensuremath{\triangleleft}

The following proposition constitutes a slight generalization of Proposition 2.2 in \cite{BdTj13b}. It establishes sufficient conditions for the existence of a switching observer exhibiting the state determination of \eqref{system:linearout}, without any a priori information concerning the initial condition. We make the following hypothesis:

\textit{Hypothesis 2.2:} There exist constants $L>1$, $c_{1},c_{2}>0$ such that \eqref{constant:c1} holds, an integer $\ell\in {\mathbb N}$, a function $g\in C^{1} ({\mathbb R}_{\ge 0};{\mathbb R})$ satisfying \eqref{map:g:properties} and a map $A:{\mathbb R}_{\ge 0}\times {\mathbb R}^{\ell}\times {\mathbb R}^{n}\times {\mathbb R}^{n}\times {\mathbb R}^{\bar{n}}\to {\mathbb R}^{n\times n}$ satisfying P1, in such a way that for every $R>0$ Hypothesis 2.1 is fulfilled, namely, both A1 and A2 hold.

\textit{Proposition 2.2:} In addition to the hypothesis of forward completeness for \eqref{system:linearout:st}, assume that system \eqref{system:linearout} satisfies Hypothesis 2.2. Then the SODP is solvable for \eqref{system:linearout} with respect to ${\mathbb R}^{n}$. \ensuremath{\triangleleft}

The proofs of Propositions 2.1 and 2.2 and the main result of Theorem 1.1 in the next section, are based on a preliminary technical result (Lemma 2.1 below) which constitutes a slight modification of the corresponding result of Lemma 2.1 in \cite{BdTj13b}. Let $\ell,m,n,\bar{n}\in {\mathbb N}$ and consider a pair $(H,A)$ of mappings:
\begin{subequations} \label{maps:HA:aux}
\begin{align}
H:=&H(t):{\mathbb R}_{\ge 0}\to {\mathbb R}^{\bar{n}\times m}; \label{map:H:aux} \\
A:=&A(t,q,x,e,y):{\mathbb R}_{\ge 0}\times {\mathbb R}^{\ell}\times {\mathbb R}^{n}\times {\mathbb R}^{m}\times {\mathbb R}^{\bar{n}}\to {\mathbb R}^{m\times m}\label{map:A:aux}
\end{align}
\end{subequations} 

\noindent where $H(\cdot)$ is continuous and $A(\cdot,\cdot,\cdot,\cdot,\cdot)$ satisfies Property P1. We make the following hypothesis:

\textit{Hypothesis 2.3:} Let $g(\cdot)\in C^{0}({\mathbb R}_{\ge 0};{\mathbb R})$ satisfying \eqref{map:g:values} and assume that for certain constant $R>0$ and for every $\xi\ge 1$, there exist a function $\beta_{R}:=\beta_{R,\xi}\in N$ and set-valued mappings $[0,\infty)\ni  t\to Y_{R}(t):=Y_{R,\xi}(t)\subset {\mathbb R}^{\bar{n}}$ and $[0,\infty)\ni t\to Q_{R}(t):=Q_{R,\xi}(t)\subset {\mathbb R}^{\ell}$ with $Y_{R}(t)\ne \emptyset$ and $Q_{R}(t)\ne \emptyset$ for all $t\ge 0$, satisfying the CP, in such a way that for every $t_{0} \ge 0$, a time-varying symmetric matrix $P_{R}:=P_{R,\xi,t_{0}}\in C^{1}([t_{0},\infty);{\mathbb R}^{m\times m})$ can be found, satisfying $P_{R}(t)\ge I_{m\times m},\forall t\ge t_{0}$ and a function $d_{R}:=d_{R,\xi,t_{0}}\in C^{0}([t_{0},\infty);{\mathbb R})$, in such a way that
\begin{align}
e'P_{R}(t)A(t,q,x,e,y)e+\tfrac{1}{2}e'\dot{P}_{R}(t)e\le-d_{R}(t)e'P_{R}(t)e,\forall t\ge t_{0},q&\in Q_{R}(t), \nonumber \\
x\in {\mathbb R}^{n},e\in\ker H(t),y\in Y_{R}(t):|x|\le\beta_{R}(t),|e|\le\xi,e'P_{R}(t)e&\ge g(t) \label{Lyapunov:ker:inequatily:aux}
\end{align}

\textit{Lemma 2.1:} Consider the pair $(H,A)$ of the time-varying mappings in \eqref{maps:HA:aux} and assume that Hypothesis 2.3 is fulfilled for certain $R>0$. Then for every $\xi\ge 1$, $t_{0} \ge 0$ and $\bar{d}_{R}:=\bar{d}_{R,\xi,t_{0}}\in C^{0}([t_{0},\infty);{\mathbb R})$ with $\bar{d}_{R}(t)<d_{R}(t)$, $\forall t\ge t_{0}$, there exists a function $\phi _{R}:=\phi_{R,\xi,t_{0}}\in C^{1}([t_{0},\infty);{\mathbb R}_{>0})$ such that
\begin{align} 
e'P_{R}(t)A(t,q,x,e,y)e+\tfrac{1}{2}e'\dot{P}_{R}(t)e\le\phi_{R}(t)|H(t)e|^{2}-\bar{d}_{R}(t)&e'P_{R} (t)e,\forall t\ge t_{0}, \nonumber \\
q\in Q_{R}(t),x\in {\mathbb R}^{n},e\in {\mathbb R}^{m},y\in Y_{R}(t):|x|\le\beta_{R}(t),|e|\le\xi,&e'P_{R}(t)e\ge g(t)\;\; \triangleleft \label{Lyapunov:inequality:aux} 
\end{align}

\section{PROOF OF THEOREM 1.1}

In this section we apply the results of Section II to prove our main result concerning the solvability of the SODP(ODP) for triangular systems \eqref{system:triangular}.

\textit{Proof of Theorem 1.1:} The proof of both statements is based on the results of Propositions 2.1 and 2.2 and is based on a generalization of the methodology employed in the proof of the main result in \cite{BdTj13b}. Without loss of generality, we may assume that, instead of Assumption H1 it holds:

\textbf{H1.'} For each $(t;x_{1},...,x_{i})\in {\mathbb R}_{\ge 0}\times {\mathbb R}^{i}$, $i=1,...,n-1$, the function $\mathbb{R}\ni z\to f_{i}(t,x_{1},...,x_{i},z)\in\mathbb{R}$ is strictly increasing.

The proof of the first statement, is based on the establishment of Hypothesis 2.2 for system \eqref{system:triangular}. Hence, we show that there exist an integer $\ell\in {\mathbb N}$, a map $A:{\mathbb R}_{\ge 0}\times {\mathbb R}^{\ell}\times {\mathbb R}^{n}\times {\mathbb R}^{n}\times {\mathbb R}\to{\mathbb R}^{n\times n}$ satisfying P1, constants $L>1$, $c_{1},c_{2} >0$ such that \eqref{constant:c1} holds and a function $g(\cdot)$ satisfying \eqref{map:g:properties}, in such a way that for each $R>0$, both A1 and A2 hold for \eqref{system:triangular}. Let $R>0$, $\xi\ge 1$ and define:
\begin{align}
F(t,x,y):=&(f_{1}(t,y,x_{2}),f_{2}(t,y,x_{2},x_{3}),...,f_{n-1}(t,y,x_{2},...,x_{n}),f_{n}(t,y,x_{2},...,x_{n}))', \nonumber \\
&(t,x,y)\in {\mathbb R}_{\ge 0} \times {\mathbb R}^{n} \times {\mathbb R} \label{rhs:var:triang}
\end{align} 

\noindent Also, for every pair of indices $(i,j)$ with $2\le j\le n$, $j-1\le i\le n$ we define the functions $\delta_{i,j}(\cdot)$ as
\begin{subequations}
\begin{align}
&\begin{array}{c} \delta_{i,2}(t,y,x_{2},...,x_{i+1},e_{2}) \\ 1\le i\le n;n\ge 2 \end{array} \nonumber \\
&\left\{\begin{array}{l}:= \frac{f_{i}(t,y,x_{2},x_{3},...,x_{i+1})-f_{i}(t,y,x_{2}-e_{2},x_{3},...,x_{i+1})}{e_{2}}, \\
\phantom{:=}\:{\rm for}\: e_{2}\ne 0 \:{\rm and}\: 2\le i\le n;n\ge 3 \\
:=\frac{f_{i}(t,y,x_{2})-f_{i}(t,y,x_{2}-e_{2})}{e_{2}},\:{\rm for}\: e_{2}\ne 0 \:{\rm and}\: i=1,2;n=2 \:{\rm or} \: i=1;n\ge 3\\ 
:=0,\:{\rm for}\: e_{2}=0 \:{\rm and}\: n,i,j \:\textup{in each case above} \end{array}\right. \nonumber \\
&(t;y;x_{2},...,x_{i+1};e_{2})\in {\mathbb R}_{\ge 0}\times\mathbb{R}\times {\mathbb R}^{i}\times\mathbb{R} \label{map:delta:i2} 
\end{align}
\begin{align}
&\begin{array}{c} \delta_{i,j}(t,y,x_{2},...,x_{i+1},e_{2},...,e_{j}) \\ 3\le j\le i\le n;n\ge 3 \end{array} \nonumber \\
&\left\{ \begin{array}{l} :=\frac{f_{i}(t,y,x_{2}-e_{2},...,x_{j-1}-e_{j-1},x_{j},x_{j+1},...,x_{i+1})-f_{i}(t,y,x_{2}-e_{2},...,x_{j-1}-e_{j-1},x_{j}-e_{j},x_{j+1},...,x_{i+1})}{e_{j}}, \\
\phantom{:=}\: {\rm for}\: e_{j} \ne 0 \:{\rm and}\: 3\le i\le n-1;3\le j\le i;n\ge 4 \:{\rm or}\: i=n;3\le j\le n-1;n\ge 4 \\
:=\frac{f_{n}(t,y,x_{2}-e_{2},...,x_{n-1}-e_{n-1},x_{n})-f_{n}(t,y,x_{2}-e_{2},...,x_{n-1}-e_{n-1},x_{n}-e_{n})}{e_{n}}, \\
\phantom{:=}\:{\rm for}\: e_{n}\ne 0 \:{\rm and}\: i=j=n;n\ge 3 \\
:=0,\:{\rm for}\: e_{j}=0, \:{\rm and}\: n,i,j \: \textup{in each case above} \end{array}\right. \nonumber \\
&(t;y;x_{2},...,x_{i+1};e_{2},...,e_{j})\in {\mathbb R}_{\ge 0}\times\mathbb{R}\times\Rat{i}\times {\mathbb R}^{j-1} \label{map:delta:ij} 
\end{align}
\begin{align}
&\begin{array}{c} \delta_{i,i+1}(t,y,x_{2},...,x_{i+1},e_{2},...,e_{i+1}) \\ 2\le i\le n-1;n\ge 3 \end{array} \nonumber \\
&\left\{\begin{array}{l} :=\frac{f_{i}(t,y,x_{2}-e_{2},...,x_{i}-e_{i},x_{i+1})-f_{i}(t,y,x_{2}-e_{2},...,x_{i}-e_{i},x_{i+1}-e_{i+1})}{e_{i+1}}, \:{\rm for}\: e_{i+1}\ne 0\\ 
:=0, \:{\rm for}\: e_{i+1}=0 \end{array}\right. \nonumber \\
&(t;y;x_{2} ,...,x_{i+1};e_{2},...,e_{i+1})\in {\mathbb R}_{\ge 0}\times\mathbb{R}\times\Rat{i}\times {\mathbb R}^{i}, \label{map:delta:iipl1}
\end{align}
\end{subequations}

\noindent (where we have used the notation $x_{i+1}|_{i=n}:=x_{n}$ and $\Rat{i}|_{i=n}:=\Rat{n-1}$, in \eqref{map:delta:i2}, \eqref{map:delta:ij}). By exploiting the  Lipschitz continuity assumption for the functions $f_{i}(\cdot)$, $i=1,...,n$ it follows that for every $2\le j\le n$, $j-1\le i\le n$,  the following properties hold:

\textbf{S1.} $\delta_{i,j}(\cdot)$ is continuous on the set
\begin{equation} \label{set:Oij}
{\mathcal O}_{i,j}:={\mathbb R}_{\ge 0}\times\mathbb{R}\times\Rat{i}\times\{(e_{2},...,e_{j})\in\Rat{j-1}:e_{j}\ne 0\}
\end{equation}

\textbf{S2.} $\delta_{i,j} (\cdot )$ is bounded on every compact subset of ${\mathbb R}_{\ge 0}\times {\mathbb R}\times\Rat{i}\times\Rat{j-1}$. (where we have used the notation $\Rat{i}|_{i=n}:=\Rat{n-1}$ in both S1 and S2)

\noindent Furthermore, from H1', \eqref{set:Oij}, \eqref{map:delta:i2} and \eqref{map:delta:iipl1} it follows:

\textbf{S3.} 
\begin{equation} \label{map:delta:iipl1:property}\\
\delta_{i,i+1}(t,y,x_{2},...,x_{i+1},e_{2},...,e_{i+1})>0, 
\forall (t,y,x_{2},...,x_{i+1},e_{2},...,e_{i+1})\in {\mathcal O}_{i,i+1},i=1,...,n-1 
\end{equation} 

\noindent From \eqref{map:delta:i2}-\eqref{map:delta:iipl1} we deduce that for $i=1,...,n$ it holds:
\begin{align}
f_{i}(t,y,x_{2},...,x_{i+1})-f_{i}(t,y,z_{2},...,z_{i+1})&=\sum_{j=2}^{i+1}\delta_{i,j}(t,y,x_{2},...,x_{i+1},x_{2}-z_{2},....,x_{j}-z_{j})(x_{j}-z_{j}) \nonumber \\
\forall t\in {\mathbb R}_{\ge 0},y&\in {\mathbb R},(x_{2} ,...,x_{i+1}),(z_{2} ,...,z_{i+1}) \in {\mathbb R}^{i} \label{rhs:difference:scalar}
\end{align}

\noindent (where we have used the notation $\sum_{j=2}^{i+1}|_{i=n}:=\sum_{j=2}^{n}$, $x_{i+1}|_{i=n}:=x_{n}$ and $z_{i+1}|_{i=n}:=z_{n}$). Also, by invoking Property S2, the following property holds. 

\textbf{S4.} There exists a function $\sigma_{R}:=\sigma_{R,\xi}\in N\cap C^{1}([0,\infty);{\mathbb R})$ satisfying:
\begin{align}
\sigma_{R}(t)\ge\sum_{i=2}^{n}\sum_{j=2}^{i}\sup\{ &|\delta_{i,j}(t,y,x_{2},...,x_{i+1},e_{2},...,e_{j})|: |y|\le\beta(t,R), \nonumber \\
&|(x_{2},...,x_{i+1} )|\le\beta(t,R),|(e_{2},...,e_{j})|\le\xi\},\forall t\ge 0 \label{map:sigma}
\end{align}

\noindent (where we have used the notation $x_{i+1}|_{i=n}:=x_{n}$). Next, consider the set-valued map $[0,\infty)\ni t\to Q_{R}(t):=Q_{R,\xi}(t)\subset {\mathbb R}^{\ell }$, $\ell:=\frac{n(n+1)}{2}$ defined as
\begin{equation} \label{map:Q:triang}
Q_{R}(t):=\{ q=(q_{1,1};q_{2,1},q_{2,2};...;q_{n,1},q_{n,2},...,q_{n,n})\in {\mathbb R}^{\ell}:|q|\le\sigma_{R}(t)\}  
\end{equation}

\noindent that obviously satisfies the CP and consider the set valued mappings $\RgeO\times (0,\xi]\ni (t,r)\to Q_{R,i}(t,r):=Q_{R,\xi,i}(t,r)\subset \mathbb{R}\times{\mathbb R}^{i}\times\Rat{i}$, $i=1,...,n-1$ given as
\begin{align}
Q_{R,i}(t,r):=\{&(y;x_{2},...,x_{i+1};e_{2},...,e_{i+1})\in  \mathbb{R}\times{\mathbb R}^{i}\times\Rat{i}:|y|\le\beta(t,R), \nonumber \\
&|(x_{2},...,x_{i+1})|\le\beta(t,R),|(e_{2},...,e_{i+1})|\le\xi,|e_{i+1}|\ge r\} \label{map:Qi:triang} 
\end{align}
\noindent 

\noindent that also satisfy both CP and the following additional property:

\textbf{S5.} For every $(t,r)\in\RgeO\times (0,\xi]$, $q\in Q_{R,i}(t,r)$ and $\varepsilon >0$, a constant $\delta>0$ can be found, such that for every $(\tilde{t},\tilde{r})\in\RgeO\times (0,\xi]$ with $|(\tilde{t},\tilde{r})-(t,r)|<\delta$, there exists $\tilde{q}\in Q_{R,i}(\tilde{t},\tilde{r})$ with $|q-\tilde{q}|<\varepsilon$.

\noindent Finally, for $i=1,...,n-1$ define:
\begin{align}
D_{R,\xi,i}(t,r):=D_{R,i}&(t,r)=\min\{\delta_{i,i+1}(t,y,x_{2},...,x_{i+1},e_{2},...,e_{i+1}): \nonumber \\
&(y;x_{2},...,x_{i+1};e_{2},...,e_{i+1})\in Q_{R,i}(t,r)\},(t,r)\in {\mathbb R}_{\ge 0}\times (0,\xi] \label{map:Delta:Ri} 
\end{align}

\noindent By exploiting \eqref{map:Delta:Ri}, Properties S1, S3, S5 and the CP property for the mappings $Q_{R,i} (\cdot ,\cdot )$, it follows that the functions $D_{R,i}(\cdot,\cdot)$, $i=1,...,n-1$ are continuous and the following hold:
\begin{align}
0<D_{R,i}(t,r)\le\delta_{i,i+1}(t,y,x_{2},...,x_{i+1},e_{2},...,e_{i+1}),\forall &(y;x_{2},...,x_{i+1};e_{2},...,e_{i+1})\in Q_{R,i}(t,r), \nonumber \\
&(t,r)\in {\mathbb R}_{\ge 0}\times (0,\xi] \label{map:Delta:Ri:positive}
\end{align} 
\begin{equation} \label{map:Delta:Ri:decreasing}
D_{R,i}(t,r_{1})\le D_{R,i}(t,r_{2}),\forall t\in {\mathbb R}_{\ge 0},r_{1},r_{2}\in (0,\xi]\: {\rm with}\: r_{1}<r_{2}
\end{equation}

\noindent Now, let $Y_{R}(\cdot)$ as given by \eqref{output:set} with
\begin{equation} \label{output:triang} 
H:=(\underbrace{1,0,...,0}_{n}) 
\end{equation} 

\noindent and notice that, due to \eqref{state:bound}, \eqref{output:set} and \eqref{output:triang}, it holds:
\begin{equation} \label{output:bound}
|y|\le\beta(t,R),\: {\rm for}\: {\rm every}\: y\in Y_{R}(t), t\ge 0
\end{equation} 

\noindent From \eqref{rhs:var:triang}, \eqref{rhs:difference:scalar}, \eqref{map:Q:triang}, \eqref{output:bound} and Property S4, it follows that for every $t\ge 0$, $y\in Y_{R} (t)$ and $x,z\in {\mathbb R}^{n} $ with $|x|\, \le \beta (t,R)$ and $|x-z|\, \le \xi $ we have:
\begin{subequations}
\begin{align}
&F(t,x,y)-F(t,z,y)=A(t,q,x,x-z,y)(x-z),\nonumber \\ 
&{\rm for}\: {\rm some}\: q\in Q_{R}(t)\: {\rm with}\: q_{i,1}=0,i=1,...,n; \label{rhs:difference:triang}
\end{align}
\end{subequations}

\noindent with 
\begin{equation}
A(t,q,x,e,y):=\left(\begin{matrix}
q_{1,1} & \delta _{1,2}(t,y,x_{2},e_{2}) & 0 & \cdots & 0 \\ 
q_{2,1} & q_{2,2} & \delta_{2,3}(t,y,x_{2},x_{3},e_{2},e_{3}) & \ddots & \vdots \\ 
\vdots & \vdots & \vdots & \ddots & 0 \\ 
q_{n-1,1} & q_{n-1,2} & q_{n-1,3} &  & \delta _{n-1,n}(t,y,x_{2},...,x_{n},e_{2},...,e_{n}) \\
q_{n,1} & q_{n,2} & q_{n,3} & \cdots & q_{n,n} \end{matrix}\right) \label{map:A:triang} 
\end{equation}

\noindent Notice that this map has the form \eqref{map:A:aux:form}, and, due to S1 and S2, satisfies Property P1. Hence, A1 is satisfied.

In order to establish A2, we prove that there exist constants $L>1$, $c_{1} ,c_{2} >0$ such that \eqref{constant:c1} holds and a function $g(\cdot )$ satisfying \eqref{map:g:properties}, in such a way that for every $R>0$, $\xi\ge 1$ and $t_{0} \ge 0$, a time-varying symmetric matrix $P_{R}:=P_{R,\xi,t_{0}} \in C^{1}([t_{0},\infty);{\mathbb R}^{n\times n})$ and a function $d_{R}:=d_{R,\xi,t_{0}}\in C^{0}([t_{0},\infty);{\mathbb R})$ can be found satisfying all conditions \eqref{map:PR:properties}, \eqref{map:dR:properties}, \eqref{Lyapunov:ker:inequality} with $H$, $A(\cdot,\cdot,\cdot,\cdot,\cdot)$, $Y_{R}(\cdot)$ and $Q_{R}(\cdot)$, as given by \eqref{output:triang}, \eqref{map:A:triang}, \eqref{output:set} and \eqref{map:Q:triang}, respectively and with $\beta(\cdot,\cdot)$ as given in \eqref{state:bound} for the case of system \eqref{system:triangular}. We proceed by induction as follows. Pick $L>1$, $c_{1}:=1$, $c_{2}:=n$ and let $g(\cdot)$ be a $C^{1}$ function satisfying \eqref{map:g:properties}. Also, let $R>0$ and for $k=2,...,n$ define:
\begin{subequations} \label{maps:triangular}
\begin{equation} \label{map:Hk}
H_{k}:=(\underbrace{1,0,...,0}_{k}),e:=(e_{n-k+1};\hat{e}')'\in {\mathbb R}\times {\mathbb R}^{k-1},\hat{e}:=(e_{n-k+2},...,e_{n})'\in {\mathbb R}^{k-1} 
\end{equation}

\noindent and consider the map $A_{k}:\RgeO\times\Rat{\ell}\times\Rat{n}\times\Rat{k}\times\mathbb{R}\to\Rat{k\times k}$ with components:
\begin{equation} \label{map:Ak}
(A_{k}(t,q,x,e,y))_{i,j} \left\lbrace\begin{array}{l}
:=q_{n-k+i,n-k+j},\:{\rm for}\: j\le i \\
:=\delta_{n-k+i,n-k+j}(t,y,x_{2},...,x_{n-k+j},0,...,0, \\
\phantom{:=}e_{n-k+1},...,e_{n-k+j}),\:{\rm for}\: j=i+1;k<n-1 \\
:=\delta_{n-k+i,n-k+j}(t,y,x_{2},...,x_{n-k+j},e_{2},...,e_{n-k+j}), \\
\phantom{:=}{\rm for}\: j=i+1;k=n-1,n \\
:=0,\:{\rm for}\: j>i+1
\end{array}\right.
\end{equation}
\end{subequations}

\noindent \textbf{Claim 1 (Induction Hypothesis):} Let $k\in {\mathbb N}$ with $2\le k\le n$. Then for $L$, $R$ and $g(\cdot)$ as above and for every $\xi\ge 1$ and $t_{0}\ge 0$, there exist a time-varying symmetric matrix $P_{R,k}:=P_{R,\xi,t_{0},k}\in C^{1}([t_{0},\infty);{\mathbb R}^{k\times k})$ and a mapping $d_{R,k}:=d_{R,\xi,t_{0},k}\in C^{0}([t_{0},\infty);{\mathbb R})$, in such a way that the following hold:
\begin{subequations}
\begin{equation}
P_{R,k}(t)>I_{k\times k},\forall t\ge t_{0}; |P_{R,k}(t_{0})|\le L; \label{map:PRk:properties}
\end{equation}
\begin{equation} \label{map:dRk:properties}
d_{R,k}(t)>n-k+1,\forall t\ge t_{0}+1;\int_{t_{1}}^{t_{2}}d_{R,k}(s)ds>-k,\forall t_{2}\ge t_{1},t_{1},t_{2}\in [t_{0},t_{0}+1] 
\end{equation}
\begin{align}
e'P_{R,k}(t)A_{k}(t,q,x,e,y)e&+\tfrac{1}{2}e'\dot{P}_{R,k}(t)e\le-d_{R,k}(t)e'P_{R,k}(t)e,\forall t\ge t_{0},q\in Q_{R}(t), \nonumber \\ 
x\in {\mathbb R}^{n},e\in\ker H_{k},y&\in Y_{R}(t):|x|\le\beta(t,R),|e|\le\xi,e'P_{R,k}(t)e\ge g(t) \label{Lyapunov:inequality:k}
\end{align}
\end{subequations}

\noindent with $H_{k}$, $A_{k}(\cdot,\cdot,\cdot,\cdot,\cdot)$, $Y_{R}(\cdot)$ and $Q_{R} (\cdot)$ as given in \eqref{map:Hk}, \eqref{map:Ak}, \eqref{output:set} and \eqref{map:Q:triang}, respectively.

By taking into account \eqref{maps:triangular}, it follows that the mappings $H_{n}$ and $A_{n}(\cdot,\cdot,\cdot,\cdot,\cdot)$ coincide with $H$ and $A(\cdot,\cdot,\cdot,\cdot,\cdot)$ as given by \eqref{output:triang} and \eqref{map:A:triang}, respectively, hence, A2 is a consequence of Claim 1 with $H:=H_{n}$ and $A(\cdot,\cdot,\cdot,\cdot,\cdot):=A_{n}(\cdot,\cdot,\cdot,\cdot,\cdot)$ and with $d_{R}:=d_{R,n}$ and $P_{R}:=P_{R,n}$ as given in  \eqref{map:dRk:properties}, \eqref{map:PRk:properties}. Indeed, relations \eqref{map:PR:properties}, \eqref{Lyapunov:ker:inequality} follow directly from \eqref{map:PRk:properties}, \eqref{Lyapunov:inequality:k} and both inequalities of \eqref{map:dR:properties} are a consequence of \eqref{map:dRk:properties} with $k=n$, if we take into account that $c_{1}=1$ and $c_{2}=n$.

\noindent \textbf{Proof of Claim 1 for $k:=2$:} For reasons of notational simplicity, we may assume that $n>3$. In that case we may define:

\begin{subequations} \label{maps:triangular:2}
\begin{equation}
H_{2}:=(1,0),e:=(e_{n-1},e_{n})'\in {\mathbb R}^{2} \label{map:H2}
\end{equation}
\begin{equation}
A_{2}(t,q,x,e,y):=\left(\begin{matrix} q_{n-1,n-1} & \delta_{n-1,n}(t,y,x_{2},...,x_{n},0,...,0,e_{n-1},e_{n}) \\ q_{n,n-1} & q_{n,n} \end{matrix}\right) \label{map:A2}
\end{equation}
\end{subequations}

\noindent Also, consider the constants $L$, $R$ and the function $g(\cdot)$ as above and let $\xi\ge 1$ and $t_{0}\ge 0$. We establish existence of a time-varying symmetric matrix $P_{R,2}:=P_{R,\xi,t_{0},2}\in C^{1} ([t_{0},\infty);{\mathbb R}^{2\times 2})$ and a mapping $d_{R,2}:=d_{R,\xi,t_{0},2}\in C^{0}([t_{0},\infty );{\mathbb R})$ in such a way that
\begin{subequations}
\begin{equation}
P_{R,2}(t)>I_{2\times 2},\forall t\ge t_{0};|P_{R,2}(t_{0})|\le L; \label{map:PR2:properties}
\end{equation}
\begin{equation}
d_{R,2}(t)>n-1,\forall t\ge t_{0}+1;\int_{t_{1}}^{t_{2}}d_{R,2}(s)ds >-2,\forall t_{2}\ge t_{1},t_{1},t_{2}\in [t_{0},t_{0} +1]; \label{map:dR2:properties}
\end{equation}
\begin{align}
e'P_{R,2}(t)A_{2}(t,q,x,e,y)e&+\tfrac{1}{2}e'\dot{P}_{R,2}(t)e\le-d_{R,2}(t)e'P_{R,2}(t)e,\forall t\ge t_{0},q\in Q_{R}(t),x\in {\mathbb R}^{n},  \nonumber \\
e=(e_{n-1},e_{n})'\in {\mathbb R}^{2},y&\in Y_{R}(t):|x|\le\beta(t,R),e\in\ker H_{2},|e|\le\xi,e'P_{R,2}(t)e\ge g(t) \label{Lyapunov:inequality:2}
\end{align}
\end{subequations}

\noindent with $H_{2}$, $A_{2}(\cdot,\cdot,\cdot,\cdot,\cdot)$, $Y_{R}(\cdot)$ and $Q_{R}(\cdot)$ as given in \eqref{map:H2}, \eqref{map:A2}, \eqref{output:set} and \eqref{map:Q:triang}, respectively. Define:
\begin{equation} \label{map:PR2} 
P_{R,2}(t):=\left(\begin{matrix} p_{R,1}(t) & p_{R}(t) \\ p_{R}(t) & L \end{matrix}\right), t\ge t_{0}
\end{equation} 

\noindent for certain $p_{R,1},p_{R}\in C^{1}([t_{0},\infty);{\mathbb R})$, to be determined in the sequel and notice, that due to \eqref{map:H2} and \eqref{map:PR2}, we have: 
\begin{equation} \label{ker:H2:equality} 
\{e\in\ker H_{2}:|e|\le\xi\:{\rm and}\:e'P_{R,2}(t)e\ge g(t)\}=\{e=(0,e_{n})':\sqrt{g(t)/L}\le |e_{n}|\le\xi\} 
\end{equation} 

\noindent Then, by taking into account \eqref{maps:triangular:2}, \eqref{map:PR2} and \eqref{ker:H2:equality}, the desired \eqref{Lyapunov:inequality:2} is written:
\begin{align}
p_{R}(t)\delta_{n-1,n}(t,y,x_{2},...,x_{n},0,...,0,e_{n})+Lq_{n,n}\le-Ld_{R,2}(t),\forall t&\ge t_{0},q\in Q_{R}(t), \nonumber \\
x\in {\mathbb R}^{n},e=(0,e_{n} )'\in {\mathbb R}^{2},y\in Y_{R}(t):|x|\le\beta(t,R),\sqrt{g(t)/L}&\le |e_{n}|\le\xi \label{Lin2:follows:from1} 
\end{align} 

\noindent By invoking \eqref{map:Q:triang}, \eqref{map:Qi:triang}, \eqref{output:bound} and the equivalence between \eqref{Lyapunov:inequality:2} and \eqref{Lin2:follows:from1}, it follows that, in order to prove \eqref{Lyapunov:inequality:2}, it suffices to determine $p_{R,1},p_{R}\in C^{1}([t_{0},\infty);{\mathbb R})$ and $d_{R,2} \in C^{0}([t_{0},\infty);{\mathbb R})$ in such a way that \eqref{map:PR2:properties} and \eqref{map:dR2:properties} are fulfilled, and further:
\begin{align} \label{Lin2:follows:from2} 
p_{R}(t)\delta_{n-1,n}(t,y,x_{2},...,x_{n},0,...,0,e_{n})&+L\sigma_{R}(t) \le-Ld_{R,2}(t),\forall t\ge t_{0}, \nonumber \\
(y;x_{2},...,x_{n};0,...,0,e_{n})&\in Q_{R,n-1}(t,\sqrt{g(t)/L})  
\end{align} 

\noindent We also require, that the candidate function $p_{R}(\cdot)$ satisfies: 
\begin{equation} \label{map:pR:keq2} 
p_{R}(t)\le 0,\forall t\ge t_{0}; p_{R}(t_{0})=0 
\end{equation} 

\noindent Then, by taking into account \eqref{map:Delta:Ri:positive}, \eqref{map:pR:keq2} and \eqref{Lin2:follows:from2}, it suffices to prove: 
\begin{equation} \label{Lin2:follows:from3} 
p_{R}(t)D_{R,n-1}(t,\sqrt{g(t)/L})+L\sigma_{R}(t)\le -Ld_{R,2}(t),\forall t\ge t_{0}  
\end{equation} 

\noindent for certain $p_{R,1},p_{R}\in C^{1}([t_{0},\infty);{\mathbb R})$ and $d_{R,2}\in C^{0}([t_{0},\infty);{\mathbb R})$ satisfying \eqref{map:PR2:properties}, \eqref{map:dR2:properties} and \eqref{map:pR:keq2}.

\noindent \textbf{Construction of $p_{R} $ and $d_{R,2} $:} First, notice that the mapping $t\to D_{R,n-1}(t,\sqrt{g(t)/L})$, $t\ge t_{0}$ is continuous, and due to \eqref{map:Delta:Ri:positive} we can find a function $\mu_{2}\in C^{1}([t_{0},\infty);{\mathbb R})$, satisfying:
\begin{equation} \label{map:mu2:propery}
0<\mu_{2}(t)\le D_{R,n-1}(t,\sqrt{g(t)/L}),\: {\rm for}\: {\rm every}\: t\ge t_{0}                                  
\end{equation}

\noindent Let
\begin{subequations} \label{constants:M2:tau2}
\begin{align}
M_{2}:=&\max\{\sigma_{R}(t):t\in [t_{0},t_{0}+\tfrac{1}{2}]\} \label{constant:M2} \\
\tau_{2}:=&\min\left\{\frac{1}{M_{2}},1\right\} \label{constant:tau2}
\end{align}
\end{subequations}

\noindent and define $\theta :=\theta_{R,\xi,t_{0}}\in C^{1}([t_{0},\infty);{\mathbb R})$, $p_{R}\in C^{1}([t_{0},\infty);{\mathbb R})$ and $d_{R,2}\in C^{0}([t_{0},\infty);{\mathbb R})$ as follows:

\noindent
\begin{align} 
\theta(t)&\left\lbrace\begin{array}{ll}
:=0,      & t=t_{0} \\
\in[0,1], & t\in\left[t_{0},t_{0}+\frac{\tau_{2}}{2}\right] \\
:=1,      & t\in\left[t_{0}+\frac{\tau_{2}}{2},\infty\right)
\end{array}\right. \label{map:theta:keq2} \\
p_{R}(t)&:=-\theta(t)\frac{L(n+\sigma_{R}(t))}{\mu_{2}(t)}, t\ge t_{0} \label{map:pR:keq2:dfn} \\
d_{R,2}(t)&\left\lbrace\begin{array}{ll}
:=-M_{2},      & t\in\left[t_{0},t_{0}+\frac{\tau_{2}}{2}\right] \\
\in[-M_{2},n], & t\in\left[t_{0}+\frac{\tau_{2}}{2},t_{0}+\tau_{2}\right] \\
:=n,           & t\in[t_{0}+\tau_{2},\infty)
\end{array}\right. \label{map:dR2}
\end{align} 

\noindent We show that \eqref{map:dR2:properties}, \eqref{map:pR:keq2} and \eqref{Lin2:follows:from3} are fulfilled, with $p_{R}(\cdot)$ and $d_{R,2}(\cdot)$ as given by \eqref{map:pR:keq2:dfn} and \eqref{map:dR2}, respectively. Indeed, \eqref{map:pR:keq2} follows directly by recalling \eqref{map:mu2:propery}, \eqref{map:theta:keq2} and \eqref{map:pR:keq2:dfn}. Both inequalities of \eqref{map:dR2:properties} are a direct consequence of \eqref{constant:tau2} and \eqref{map:dR2}. We next show that \eqref{Lin2:follows:from3} holds as well, with $p_{R}(\cdot)$ and $d_{R,2}(\cdot)$ as above and consider two cases:

\noindent \textbf{Case 1:} $t\in [t_{0},t_{0}+\tfrac{\tau_{2}}{2}] $. In that case, \eqref{Lin2:follows:from3} follows directly from \eqref{map:pR:keq2}, \eqref{map:mu2:propery}, \eqref{constants:M2:tau2} and \eqref{map:dR2}.

\noindent \textbf{Case 2:} $t\in [t_{0}+\tfrac{\tau_{2}}{2},\infty)$. Then from \eqref{map:mu2:propery}, \eqref{map:theta:keq2}, \eqref{map:pR:keq2:dfn} and \eqref{map:dR2} it follows that:
\begin{equation*}
-\frac{L(n+\sigma_{R}(t))}{\mu_{2}(t)}D_{R,n-1}(t,\sqrt{g(t)/L})+L\sigma_{R}(t)\le-Ln\le -Ld_{R,2}(t)
\end{equation*}

\noindent namely, \eqref{Lin2:follows:from3} again holds for all $t\in [t_{0}+\tfrac{\tau _{2}}{2},\infty)$.

\noindent We therefore conclude that \eqref{Lin2:follows:from3} is fulfilled for all $t\ge t_{0}$.

\noindent Finally, the construction of $p_{R,1}(\cdot)$ included in \eqref{map:PR2}, is the same with that given in proof of Theorem 1.1 in \cite{BdTj13b} and is omitted. This completes the proof of Claim 1 for $k=2$.

\noindent \textbf{Proof of Claim 1 (general step of induction procedure):} Assume  now that Claim 1 is fulfilled for certain integer $k$ with $2\le k<n$. We prove that Claim 1 also holds for $k:=k+1$. Consider the pair $(H,A)$ as given in \eqref{maps:HA:aux} with $H(t):=H_{k}$, $A(t,q,x,e,y):=A_{k}(t,q,x,e,y)$, $\ell =\frac{n(n+1)}{2}$, $m:=k$, $n:=n$ and $\bar{n}:=1$, where $H_{k}$ and $A_{k}$ are defined by \eqref{map:Hk} and \eqref{map:Ak}, respectively. Notice, that the map $A_{k}$ as given by \eqref{map:Ak} has the form \eqref{map:A:aux:form} and due to S1 and S2, satisfies Property P1. Hence, by the first inequality of \eqref{map:PRk:properties} and \eqref{Lyapunov:inequality:k}, we conclude that Hypothesis 2.3 of the previous section holds, with $R$ and $g(\cdot)$ as above, $Y_{R}(\cdot)$, $Q_{R}(\cdot)$ and $\beta_{R}(\cdot):=\beta(\cdot,R)$ as given in \eqref{output:set}, \eqref{map:Q:triang} and \eqref{state:bound}, respectively, and with $d_{R}(\cdot):=d_{R,k}(\cdot)$ and $P_{R}(\cdot):=P_{R,k}(\cdot)$ as given in \eqref{map:PRk:properties}, \eqref{map:dRk:properties}. Finally, for every $\xi\ge 1$ and $t_{0}\ge 0$, consider the function $\bar{d}_{R,k}:=\bar{d}_{R,\xi,t_{0},k}$ defined as:
\begin{equation} \label{map:dRbk} 
\bar{d}_{R,k}(t):=d_{R,k}(t)-\tfrac{1}{2},t\ge t_{0}  
\end{equation} 

\noindent which satisfies $\bar{d}_{R,k}(t)<d_{R,k}(t)$ for all $t\ge t_{0} $. It follows that all requirements of Lemma 2.1 are fulfilled and therefore, there exists a function $\phi_{R,k}:=\phi_{R,\xi,t_{0},k}\in C^{1}([t_{0},\infty);{\mathbb R}_{>0})$ such that
\begin{align}
e'P_{R,k}(t)A_{k}(t,q,x,e,y)e+\tfrac{1}{2}e'\dot{P}_{R,k}(t)&e\le\phi_{R,k}(t)|H_{k} e|^{2}-\bar{d}_{R,k}(t)e'P_{R,k} (t)e,\forall t\ge t_{0}, \nonumber \\ 
q\in Q_{R}(t),x\in {\mathbb R}^{n},e\in {\mathbb R}^{k},y\in Y_{R}(t)&:|x|\le\beta(t,R),|e|\le\xi,e'P_{R,k}(t)e\ge g(t) \label{Lyapunov:inequality:aux:k} 
\end{align}

\noindent Furthermore, due to \eqref{map:dRk:properties} and \eqref{map:dRbk}, the map $\bar{d}_{R,k}(\cdot)$ satisfies:
\begin{equation} \label{map:dRbk:properties}
\bar{d}_{R,k}(t)>n-k+\tfrac{1}{2},\forall t\ge t_{0}+1;\int_{t_{1}}^{t_{2}}\bar{d}_{R,k} (s)ds >-(k+\tfrac{1}{2}),\forall t_{2}\ge t_{1},t_{1},t_{2}\in [t_{0},t_{0}+1] 
\end{equation}

In the sequel, we exploit \eqref{Lyapunov:inequality:aux:k} and \eqref{map:dRbk:properties}, in order to establish that Claim 1 is fulfilled for $k=k+1$. Specifically, for the same $L$, $R$ and $g(\cdot)$ as above and for any $\xi\ge 1$ and $t_{0}\ge 0$, we show that there exist a time-varying symmetric matrix $P_{R,k+1}\in C^{1}([t_{0},\infty);{\mathbb R}^{(k+1)\times (k+1)})$ and a map $d_{R,k+1}\in C^{0}([t_{0},\infty);{\mathbb R})$, such that both \eqref{map:PRk:properties} and \eqref{map:dRk:properties} are fulfilled with $k=k+1$ and further:
\begin{align}
e'P_{R,k+1}(t)A_{k+1}(t,q,x,e,y)e+\tfrac{1}{2}e'\dot{P}_{R,k+1}(t)e&\le -d_{R,k+1}(t)e'P_{R,k+1}(t)e,\forall t\ge t_{0},q\in Q_{R}(t), \nonumber \\
x\in {\mathbb R}^{n} e\in\ker H_{k+1},y\in Y_{R}(t):|x|&\le\beta(t,R),|e|\le\xi, e'P_{R,k+1}(t)e\ge g(t) \label{Lyapunov:inequality:kpl1} 
\end{align}

\noindent where
\begin{subequations} \label{maps:triangular:kpl1}
\begin{equation}
H_{k+1}:=(\underbrace{1,0,...,0}_{k+1}),e:=(e_{n-k};\hat{e}')'\in {\mathbb R}\times {\mathbb R}^{k},\hat{e}:=(e_{n-k+1},...,e_{n})'\in {\mathbb R}^{k} \label{map:Hkpl1} 
\end{equation}

\noindent the components of the map $A_{k+1}:\RgeO\times\Rat{\ell}\times\Rat{n}\times\Rat{k+1}\times\mathbb{R}\to\Rat{(k+1)\times (k+1)}$ are defined as:
\begin{equation} \label{map:Akpl1}
(A_{k+1}(t,q,x,e,y))_{i,j}\left\lbrace\begin{array}{l}
:=q_{n-k-1+i,n-k-1+j},\:{\rm for}\: j\le i \\
:=\delta_{n-k-1+i,n-k-1+j}(t,y,x_{2},...,x_{n-k+i},0,...,0, \\
\phantom{:=}e_{n-k},...,e_{n-k+i}),\:{\rm for}\: j=i+1;k<n-2 \\
:=\delta_{n-k-1+i,n-k-1+j}(t,y,x_{2},...,x_{n-k+i},e_{2},...,e_{n-k+i}), \\
\phantom{:=}\:{\rm for}\: j=i+1;k=n-2,n-1 \\
:=0,\:{\rm for}\: j>i+1
\end{array}\right.
\end{equation}
\begin{equation}
P_{R,k+1}(t):=\left(\begin{matrix} p_{R,1}(t) & \begin{matrix} p_{R}(t) & 0 & \cdots & 0 \end{matrix} \\ \begin{matrix} p_{R}(t) \\ 0 \\ \vdots \\ 0 \end{matrix} & \boxed{P_{R,k}(t)} \end{matrix}\right) \label{map:PRkpl1}
\end{equation}
\end{subequations}

\noindent and where $Y_{R}(\cdot)$, $Q_{R}(\cdot)$ are given in \eqref{output:set} and \eqref{map:Q:triang}, respectively. Again, for reasons of notational simplicity, we may assume that $k<n-1$. It then follows from \eqref{map:Ak} and \eqref{map:Akpl1} that for every $e=(0,\hat{e}')'=(0,e_{n-k+1},...,e_{n})'\in \ker H_{k+1}:=(\underbrace{1,0,...,0}_{k+1})$, the map $A_{k+1}(\cdot,\cdot,\cdot,\cdot,\cdot)$ takes the form:

\begin{equation} \label{map:Aklp1:at:kernel}
A_{k+1}(t,q,x,e,y)=\left(\begin{matrix} q_{n-k,n-k} & \begin{matrix} \delta_{n-k,n-k+1}(t,y,x_{2},...,x_{n-k+1},0,...,0,e_{n-k+1}) & 0 & \cdots & 0 \end{matrix} \\ \begin{matrix} q_{n-k+1,n-k} \\ \vdots \\ q_{n,n-k} \end{matrix} & \boxed{A_{k}(t,q,x,\hat{e},y)} \end{matrix}\right) 
\end{equation}

\noindent Let $\xi\ge 1$ and $t_{0}\ge 0$. We determine functions $p_{R,1},p_{R}\in C^{1}([t_{0},\infty);{\mathbb R})$ and $d_{R,k+1} \in C^{0} ([t_{0} ,\infty );{\mathbb R})$ such that \eqref{map:PRk:properties} and \eqref{map:dRk:properties} are fulfilled with $k=k+1$, and further \eqref{Lyapunov:inequality:kpl1} holds, with $H_{k+1} $, $A_{k+1} (\cdot ,\cdot ,\cdot ,\cdot ,\cdot )$ and $P_{R,k+1} (\cdot )$ as given by \eqref{maps:triangular:kpl1}. By taking into account \eqref{map:Aklp1:at:kernel}, it follows that \eqref{Lyapunov:inequality:kpl1} is equivalent to:
\begin{align}
e_{n-k+1}^{2}p_{R}(t)&\delta_{n-k,n-k+1}(t,y,x_{2},...,x_{n-k+1},0,...,0,e_{n-k+1})+\hat{e}'P_{R,k}(t)A_{k}(t,q,x,\hat{e},y)\hat{e} \nonumber \\
+\tfrac{1}{2}\hat{e}'\dot{P}_{R,k}(t)&\hat{e}\le -d_{R,k+1}\hat{e}'P_{R,k}(t)\hat{e},\forall t\ge t_{0},q\in Q_{R}(t),x\in {\mathbb R}^{n},e:=(e_{n-k};\hat{e}')'\in {\mathbb R}\times {\mathbb R}^{k}, \nonumber \\
y\in Y_{R} (t)&:|x|\le\beta(t,R),e\in\ker H_{k+1},|e|\le\xi,e'P_{R,k+1}(t)e\ge g(t)  \label{Linkpl1:follows:from1}
\end{align} 

\noindent Notice that, according to \eqref{map:Hkpl1} and \eqref{map:PRkpl1}, we have $e'P_{R,k+1}(t)e=\hat{e}'P_{R,k}(t)\hat{e}$ for every $e=(0,\hat{e}')'=(0,e_{n-k+1},...,e_{n})'\in \ker H_{k+1}$, thus, by taking into account  \eqref{map:Hk} and \eqref{Lyapunov:inequality:aux:k}, it suffices, instead of \eqref{Linkpl1:follows:from1}, to show that
\begin{align}
e_{n-k+1}^{2}&(p_{R}(t)\delta_{n-k,n-k+1}(t,y,x_{2},...,x_{n-k+1},0,...,0,e_{n-k+1})+\phi_{R,k}(t)) \nonumber \\
\le &(\bar{d}_{R,k}(t)-d_{R,k+1}(t))\hat{e}'P_{R,k}(t)\hat{e},\forall t\ge t_{0},x\in {\mathbb R}^{n},\hat{e}\in {\mathbb R}^{k}, \nonumber \\
&y\in Y_{R}(t):|x|\le\beta(t,R),|\hat{e}|\le\xi,\hat{e}'P_{R,k}(t)\hat{e}\ge g(t) \label{Linkpl1:follows:from2}
\end{align}

\noindent \textbf{Establishment of \eqref{Linkpl1:follows:from2} plus \eqref{map:PRk:properties} and \eqref{map:dRk:properties} for $k=k+1$:} We impose the following additional requirements for the candidate functions $p_{R}(\cdot)$ and $d_{R,k+1}(\cdot)$:
\begin{subequations} \label{maps:pR:dRkpl1}
\begin{align}
p_{R}(t)&\le 0,\forall t\ge t_{0};p_{R}(t_{0})=0; \label{map:pR:keqkpl1prop} \\
d_{R,k+1}(t)&\le\bar{d}_{R,k}(t),\forall t\ge t_{0} \label{dRkpl1:lt:dRbk}
\end{align}
\end{subequations}

\noindent Then, by taking into account \eqref{dRkpl1:lt:dRbk} and the fact that the desired inequality in \eqref{Linkpl1:follows:from2} should be valid for those $\hat{e}\in {\mathbb R}^{k}$ for which $|\hat{e}|\le\xi$ and $\hat{e}'P_{R,k}(t)\hat{e}\ge g(t)$, it follows that, in order to show \eqref{Linkpl1:follows:from2} and that \eqref{map:PRk:properties}, \eqref{map:dRk:properties} are valid with $k=k+1$, it suffices to show that
\begin{align}
&e_{n-k+1}^{2}(p_{R}(t)\delta_{n-k,n-k+1}(t,y,x_{2},...,x_{n-k+1},0,...,0,e_{n-k+1})+\phi_{R,k}(t))\nonumber \\
\le &(\bar{d}_{R,k}(t)-d_{R,k+1}(t))g(t), \forall t\ge t_{0},x\in\Rat{n},e_{n-k+1}\in {\mathbb R},y\in Y_{R}(t):|x|\le\beta(t,R),|e_{n-k+1}|\le\xi \label{Linkpl1:follows:from3}
\end{align}

\noindent for suitable functions $p_{R,1},p_{R}\in C^{1}([t_{0},\infty);{\mathbb R})$ and $d_{R,k+1}\in C^{0}([t_{0},\infty);{\mathbb R})$, in such a way that \eqref{map:PRk:properties}, \eqref{map:dRk:properties} hold with $k=k+1$, and in addition $p_{R}(\cdot)$ and $d_{R,k+1}(\cdot)$ satisfy \eqref{maps:pR:dRkpl1}. We proceed to the explicit construction of these functions. Due to \eqref{map:Qi:triang}, \eqref{map:Delta:Ri:positive}, \eqref{output:bound}, \eqref{map:pR:keqkpl1prop} and the fact that, due to requirement \eqref{dRkpl1:lt:dRbk}, equation \eqref{Linkpl1:follows:from3} holds trivially for $e_{n-k+1}=0$, it suffices, instead of \eqref{Linkpl1:follows:from3}, to show that
\begin{align} 
r^{2}&(p_{R}(t)D_{R,n-k}(t,r)+\phi_{R,k}(t)) \nonumber \\
\le&(\bar{d}_{R,k}(t)-d_{R,k+1}(t))g(t),\forall t\ge t_{0},r\in (0,\xi] \label{Linkpl1:follows:from4}
\end{align} 

\noindent \textbf{Construction of the mappings $p_{R}$ and $d_{R,k+1}$:} Let 
\begin{subequations} \label{constants:Mkpl1:taukpl1}
\begin{align}
M_{k+1}:=&\max\left\{|\bar{d}_{R,k}(t)|+\frac{1}{4}+\frac{\xi^{2}\phi_{R,k}(t)}{g(t)}:t\in [t_{0},t_{0}+{\tfrac{1}{2}}]\right\} \label{constant:Mkpl1} \\
\tau_{k+1}:=&\min\left\{\frac{1}{4M_{k+1}},\frac{1}{2}\right\} \label{constant:taukpl1}
\end{align}
\end{subequations}

\noindent and consider a function $\theta:=\theta_{R,\xi,t_{0}}\in C^{1}([t_{0},\infty );{\mathbb R})$ defined as:

\noindent 
\begin{equation} \label{map:theta:keqkpl1} 
\theta (t)\left\{\begin{array}{ll} :=0, & t=t_{0} \\ \in [0,1], & t\in [t_{0},t_{0}+\tfrac{\tau_{k+1}}{2}] \\ :=1, & t\in [t_{0}+\tfrac{\tau_{k+1}}{2},\infty) \end{array}\right.  
\end{equation} 

\noindent By taking into account \eqref{constants:Mkpl1:taukpl1}, it follows that:

\noindent 
\begin{equation} \label{map:dRbk:property2} 
\bar{d}_{R,k}(t)-\tfrac{1}{4}\ge -M_{k+1},\forall t\in [t_{0},t_{0}+\tau_{k+1}] 
\end{equation}

\noindent hence, by exploiting \eqref{map:dRbk:property2}, we can construct a function $d_{R,k+1}\in C^{0}([t_{0},\infty);{\mathbb R})$, satisfying:
\begin{align} \label{map:dRkpl1:dfn} 
&d_{R,k+1}(t) \nonumber \\
&\left\{\begin{array}{ll} :=-M_{k+1}, & t\in [t_{0},t_{0}+\tfrac{\tau_{k+1}}{2}] \\ \in [-M_{k+1},\bar{d}_{R,k}(t)-\tfrac{1}{4}], & t\in [t_{0}+\tfrac{\tau_{k+1}}{2},t_{0}+\tau _{k+1}] \\ :=\bar{d}_{R,k}(t)-\tfrac{1}{4}, & t\in [t_{0}+\tau_{k+1},\infty ) \end{array}\right.  
\end{align}
 
\noindent Notice that \eqref{dRkpl1:lt:dRbk}, follows from \eqref{map:dRbk:property2} and \eqref{map:dRkpl1:dfn}. Also, define:
\begin{equation} \label{map:zeta:dfn} 
\zeta(t):=\frac{1}{2}\sqrt{\frac{g(t)}{\phi_{R,k}(t)}}, t\ge t_{0}  
\end{equation} 

\noindent Due to \eqref{map:Delta:Ri} and \eqref{map:Delta:Ri:positive}, the map $t\to D_{R,n-k}(t,\zeta(t))$, $t\ge t_{0}$ is continuous and there exists a function $\mu_{k+1}\in C^{1}([t_{0},\infty);{\mathbb R})$ satisfying:
\begin{equation}\label{map:mukpl1:property}
0<\mu_{k+1}(t)\le D_{R,n-k}(t,\zeta(t)),\:{\rm for}\: {\rm every}\: t\ge t_{0}
\end{equation}

\noindent Finally, define $p_{R}\in C^{1}([t_{0},\infty);{\mathbb R})$ as:
\begin{equation} \label{map:pR:keqkpl1} 
p_{R}(t):=-\frac{\theta(t)\phi_{R,k}(t)}{\mu_{k+1}(t)}, t\ge t_{0}  
\end{equation} 

\noindent which due to \eqref{map:theta:keqkpl1} and \eqref{map:mukpl1:property}, satisfies \eqref{map:pR:keqkpl1prop}.

\noindent \textbf{Proof of \eqref{Linkpl1:follows:from4}: }We consider two cases:

\noindent \textbf{Case 1:} $t\in [t_{0},t_{0}+\tfrac{\tau_{k+1}}{2}]$. By taking into account \eqref{constant:Mkpl1}, it follows that $M_{k+1}\ge -\bar{d}_{R,k}(t)+\frac{\phi_{R,k}(t)\xi^{2}}{g(t)}$ for every $t\in [t_{0},t_{0}+\tfrac{1}{2}]$, which in conjunction with \eqref{constant:taukpl1} and \eqref{map:dRkpl1:dfn} imply:
\begin{equation} \label{map:dRkpl1:ineq} 
\bar{d}_{R,k}(t)-d_{R,k+1}(t)\ge\frac{\phi_{R,k}(t)\xi^{2}}{g(t)},\forall t\in \left[t_{0},t_{0}+\frac{\tau_{k+1}}{2}\right] 
\end{equation}

\noindent Hence, from \eqref{map:Delta:Ri:positive}, \eqref{map:pR:keqkpl1prop} and \eqref{map:dRkpl1:ineq} we deduce that
\begin{align*}
&r^{2}(p_{R}(t)D_{R,n-k}(t,r)+\phi_{R,k}(t))\le r^{2}\phi_{R,k}(t) \\
\le&\frac{\xi^{2}\phi_{R,k}(t)}{g(t)}g(t)\le (\bar{d}_{R,k}(t)-d_{R,k+1}(t))g(t),\forall r\in (0,\xi]
\end{align*}

\noindent which implies \eqref{Linkpl1:follows:from4} for $t\in [t_{0},t_{0}+\tfrac{\tau_{k+1}}{2}]$.

\noindent \textbf{Case 2:} $t\in [t_{0}+\tfrac{\tau_{k+1}}{2},\infty)$. We consider two further subcases.

\noindent \textbf{Subcase 1: $0<r\le\zeta(t)$}. Due to \eqref{map:dRkpl1:dfn}, it holds $\bar{d}_{R,k}(t)-d_{R,k+1}(t)\ge \tfrac{1}{4}$ for every $t\in [t_{0}+\tfrac{\tau_{k+1}}{2},\infty)$, hence, by exploiting \eqref{map:Delta:Ri:positive}, \eqref{map:pR:keqkpl1prop} and \eqref{map:zeta:dfn} we have
\begin{equation} \label{Linkpl1:follows:from4:ok} 
r^{2}(p_{R}(t)D_{R,n-k}(t,r)+\phi _{R,k}(t))\le\zeta^{2}(t)\phi_{R,k}(t)\le(\bar{d}_{R,k}(t)-d_{R,k+1}(t))g(t) 
\end{equation}

\noindent \textbf{Subcase 2: $\zeta(t)\le r\le\xi$}. By taking into account \eqref{dRkpl1:lt:dRbk}, \eqref{map:theta:keqkpl1}, \eqref{map:Delta:Ri:decreasing}, \eqref{map:zeta:dfn}, \eqref{map:mukpl1:property} and \eqref{map:pR:keqkpl1}, we deduce that
\begin{align*}
&r^{2}(p_{R}(t)D_{R,n-k}(t,r)+\phi_{R,k}(t))\le r^{2}\phi_{R,k}(t) \\ 
\times &\left(-\frac{D_{R,n-k}(t,\zeta(t))}{\mu_{k+1}(t)}+1\right)\le 0\le (\bar{d}_{R,k}(t)-d_{R,k+1}(t))g(t)
\end{align*}

\noindent The latter, in conjunction with \eqref{Linkpl1:follows:from4:ok} asserts that \eqref{Linkpl1:follows:from4} is fulfilled for every $t\in [t_{0}+\tfrac{\tau_{k+1}}{2},\infty)$. Both cases above guarantee that \eqref{Linkpl1:follows:from4} holds for all $t\in [t_{0},\infty)$ as required.

\noindent \textbf{Proof of \eqref{map:PRk:properties} and \eqref{map:dRk:properties} for $k=k+1$:}  The proof of \eqref{map:dRk:properties} is the same with that given in proof of Theorem 1.1 in \cite{BdTj13b} and is omitted. Finally, the proof of \eqref{map:PRk:properties} is based on the construction of the map $p_{R,1}(\cdot)$ as involved in \eqref{map:PRkpl1}, and is also the same with that given in proof of Theorem 1.1 in \cite{BdTj13b}.

We have shown that all requirements of Claim 1 hold, which, as was pointed out, establishes that for every $R>0$ Hypothesis A2 is fulfilled. We therefore conclude, that for system \eqref{system:triangular} Hypothesis 2.2 is satisfied hence, by invoking the result of Proposition 2.2, it follows that the SODP is solvable for \eqref{system:triangular} with respect to \textit{${\mathbb R}^{n} $.} The establishment of the second statement of Theorem 1.1, follows directly from Proposition 2.1. \ensuremath{\Box}

\textit{Example:} As an illustrative example of Theorem 1.1, consider the two dimensional polynomial system:
\begin{subequations} \label{system:example}
\begin{align}
\dot{x}_{1}&=x_{1}-x_{1}^{3}+x_{1}^{2}x_{2}+\frac{3}{2}x_{1}x_{2}^{2}+x_{2}^{3} \nonumber \\
\dot{x}_{2}&=-x_{1}^{3}-x_{1}x_{2}^{2}+x_{2}-x_{2}^{3} \\
y&=x_{1}
\end{align}
\end{subequations}

\noindent It is easy to check that system \eqref{system:example} satisfies all conditions of Theorem 1.1, therefore the SODP is solvable for \eqref{system:example} with respect to $\Rat{2}$.

\end{document}